\documentclass[reqno,12pt]{article}
\usepackage{amsmath}
\usepackage{amssymb}
\usepackage{latexsym}
\usepackage{amsfonts}
\usepackage{amsthm}

\newcommand{\0}{{\bf 0}}
\newcommand{\R}{\mathbb{R}}
\newcommand{\indicator}{{\mathbf 1}}
\newcommand{\Z}{\mathbb{Z}}

\newcommand{\z}{\mathbb{Z}}
\newcommand{\N}{\mathbb{N}}
\newcommand{\Exp}{{\bf E}}
\newcommand{\good}{{ V}}

\newcommand{\dens}{\mathrm{density\,}}

\newcommand{\ball}{R}
\newcommand{\sides}{\text{Sides$({\fpp})$}}

\newcommand{\E}{\mbox{$\bf E$}}
\newcommand{\booze}{Busemann }

\newcommand{\zd}{{\mathbb Z}^d}

\renewcommand{\P}{{\bf P}}

\newcommand{\ez}{{\bf 1}}

\newcommand{\olle}{{H\"{a}ggstr\"{o}m }}
\newcommand{\fpp}{{\mu}}
\newcommand{\Ball}{{\text{Ball}}}
\newcommand{\mug}{\mbox{C}}
\newcommand{\edges}{\mbox{Edges}}

\newcommand{\be}{\begin{equation}}
\newcommand{\ee}{\end{equation}}
\newcommand{\prob}{{\bf P}}

\newtheorem{thm}{Theorem}[section]
\newtheorem{theorem}[thm]{Theorem}

\newtheorem{lemma}[thm]{Lemma}
\newtheorem{corollary}[thm]{Corollary}
\newtheorem{cor}[thm]{Corollary}

 \newenvironment{pfof}[1]
{\par\vskip2\parsep\noindent{\sc Proof of\ #1. }}{{\hfill $\Box$}
\par\vskip2\parsep}

\begin{document}

\title{Geodesics in First Passage Percolation}
\author{Christopher Hoffman}
\date{}
\maketitle
\renewcommand{\thefootnote}{}
\renewcommand{\thefootnote}{\arabic{footnote}}

\abstract{ We consider a wide class of ergodic first passage
percolation processes on $\Z^2$ and prove that there exist at
least four one-sided geodesics a.s.  We also show that coexistence
is possible with positive probability in a four color Richardson's
growth model.  This improves earlier results of \olle and Pemantle
\cite{HP}, Garet and Marchand \cite {G} and Hoffman \cite{hoffman}
who proved that first passage percolation has at least two
geodesics and that coexistence is possible in a two color
Richardson's growth model.}
\section{Introduction}

\subsection{First passage percolation}
First passage percolation is a process introduced by Hammersley
and Welsh  as a time dependent model for the passage of a fluid
through a porous medium which has provided a large number of
problems of probabilistic interest with excellent physical
motivation \cite{HW}. Study of this model led to the development
of the ergodic theory of subadditive processes by Kingman
\cite{king}. It also has links to mathematical biology through
Richardson's growth model \cite{HP}. A good overview of first
passage percolation is contained in \cite{K}.

Let $\fpp$ be a stationary measure on $[0,\infty)^{\edges(\zd)}$
and let $\omega$ be a realization of $\fpp$. For any $x$ and $y$
we define $\tau(x,y)$, the {\bf passage time from $x$ to $y$}, by
$$\tau(x,y)=\inf\sum{\omega(v_{i},v_{i+1})}$$
where the sum is taken over all of the edges in the path and the
$\inf$ is taken over all paths connecting $x$ to $y$.  The time
minimizing path from $x$ to $y$ is called a {\bf geodesic}.

An infinite path $v_1, v_2, \ldots$ is called a {\bf geodesic} if
for all $0<i<j$
\[
\tau(v_i,v_j)=\sum_{k=i}^{j-1}\omega(v_k, v_{k+1}).
\]
 In
this paper we prove that for a very general class of first passage
percolation processes that there exist at least four disjoint
infinite geodesics a.s.

For notational reasons it will often be convenient to think of
$\tau$ as a function defined on $\R^2 \times \R^2$ by setting
$$\tau(x+u,y+v)=\tau(x,y)$$
for any $x,y \in \Z^2$ and any $u,v \in [-\frac12,\frac12)^2.$ For
any $x \in \R^2$ and $S\subset R^2$ we write
$$\tau(x,S)=\inf_{y \in S}\tau(x,y).$$
The most
basic result from first passage percolation is the shape theorem.
Define
$$\ball(t)=\{v:\ \tau(\0,v)\leq t\}.$$
The shape theorem says that there is a nonempty set $\ball$ such
that (modulo the boundary) $\frac{ \ball(t)}{t}$ converges to
$\ball$ a.s.
\begin{thm} \cite{B} \label{shape1}
Let $\fpp$ be stationary and ergodic, where the distribution on
any edge has finite $d+\epsilon$ moment for $\epsilon>0$. There
exists a closed set $\ball$ which is nonempty, convex, and
symmetric about reflection through the coordinate axis such that
for every $\epsilon>0$
$$P\left(\exists \ T:\ (1-\epsilon)\ball \subset \frac{ \ball(t)}{t}\subset (1+\epsilon)\ball\
    \text{ for all $t>T$}\right)=1.$$
\end{thm}
This theorem is an example of a subadditive ergodic theorem. In
general, little is known about the shape of $\ball$ other than it
is convex and symmetric. Cox and Durrett have shown that there are
nontrivial product measures such that the boundary of $\ball$
contains a flat piece yet it is neither a square nor a diamond
\cite{DL}. However for any nonempty, convex, and symmetric set
$\ball$ there exist a stationary measure $\fpp$ such that the
shape for $\fpp$ is $\ball$ \cite{H}.


Another widely studied aspect of first passage percolation are
geodesics.  We let $G(x,y)$ be the geodesic connecting $x$ and
$y$. Define

$$\Gamma(x)=\cup_{y \in \zd}
    \{e \in G(x,y)\}.$$
We refer to this as the {\bf tree of infection of $x$}. We define
$K(\Gamma(x))$ to be the number of topological ends in
$\Gamma(x)$.

Newman has conjectured that for a large class of $\fpp$,
$|K(\Gamma(\0))|=\infty$ a.s. \cite{N}  \olle and Pemantle proved
that if $d=2$, $\fpp$ is i.i.d.\ and $\omega(e)$ has exponential
distribution then with positive probability $|K(\Gamma(\0))|>1$.
In independent work Garet and Marchand \cite{G} and Hoffman
\cite{hoffman} extended this result in two directions. Their
results apply to a wide class of ergodic measures $\fpp$ on any
$d\geq2$.

Newman has proved that if $\fpp$ is i.i.d.\ and $\ball$ has
certain properties then $|K(\Gamma(\0))|=\infty$ a.s. \cite{N}.
Although these conditions are plausible there are no known
measures $\fpp$ with $S$ that satisfy these conditions.  In this
paper we prove an analogous theorem but with a much weaker
condition on $\ball$. Unfortunately even this weaker condition,
that $\partial \ball$ is not a polygon, hasn't been verified for
any version of i.i.d.\ first passage percolation.

Now we will introduce some more notation which will let us list
the conditions that we place on $\fpp$ for the rest of this paper.
We say that $\fpp$ has {\bf unique passage times} for all $x$ and
$y\neq z$
$$\prob(\tau(x,y)\neq \tau(x,z))=1.$$

Now we are ready to define the class of measures that we will work
with. We say that $\fpp$ is {\bf good} if
\begin{enumerate}
\item $\fpp$ is ergodic,
\item $\fpp$ has all the symmetries of $\zd$, \label{cont0}
\item $\fpp$ has unique passage times,
\item the distribution of $\fpp$ on any edge has finite
$2+\epsilon$ moment for some $\epsilon>0$
\item $\ball$ is bounded.
\end{enumerate}
Throughout the rest of the paper we will assume that $\mu$ is
good. Unfortunately there is no general necessary and sufficient
condition to determine when the shape $\ball$ is bounded and therefore there is no general condition for $\mu$ to be good. See
\cite{H} for examples.  However if $\mu$ is i.i.d\ and the distribution on any edge is continuous with finite $2+\epsilon$ moment then $\mu$ is good.  See Theorem 4.3 in \cite{gratuitous} for more information about conditions that imply $\mu$ is good in the case that $\mu$ is stationary but not i.i.d.

\subsection{Spatial Growth Models}
Richardson's growth model, a simple competition model between
diseases, was introduced by \olle and Pemantle \cite{HP}.  The
rules for this model are as follows.  Each vertex $z \in \z^2$ at
each time $t\geq 0$ is either infected by one of $k$ diseases
($z_t\in \{1,\dots,k\}$) or is uninfected ($z_t=0$). Initially for
each disease there is one vertex which is infected by that
disease. All other vertices are initially uninfected.
 Once a vertex is infected one of the  diseases it
stays infected by that disease for all time and is not infected by
any disease. All of the diseases spread from sites they have
already infected to neighboring uninfected sites at some rate.

We now explain the relationship between first passage percolation
and Richardson's growth models. For any $\omega \in
[0,\infty)^{\edges(\zd)}$ with unique passage times and any
$x_1,\dots,x_k \in \zd$ we can project $\omega$ to
$\tilde{\omega}_{x_1,\dots,x_k} \in
\left(\{0,1,\dots,k\}^{\zd}\right)^{[0,\infty)}$ by
$$\tilde{\omega}_{x_1,\dots,x_k}(z,t)=\left\{%
\begin{array}{ll}
     i& \hbox{if $\tau(x_i,z)\leq t$ and $\tau(x_i,z) <\tau (x_j,z)$ for all $i \neq j$}; \\
     0& \hbox{else.} \\
\end{array}%
\right.
$$
If $\fpp$ has unique passage times then $\fpp$ projects onto a
measure on
 $\left(\{0,1, \dots, k\}^{\zd}\right)^{[0,\infty)}$.
It is clear that the models start with a single vertex in states 1
through $k$.
Vertices in states $i>0$ remain in their states forever, while
vertices in state 0 which are adjacent to a vertex in state i can
switch to state i. We think of the vertices in states $i>0$ as
infected with one of $k$ infections while the vertices in state 0
are considered uninfected.

For this model it is most common to choose $\fpp$ to be i.i.d.\
with an exponential distribution on each edge.  This makes the
spatial growth process Markovian.

As each $z \in \zd$ eventually changes to some state $i>0$ and
then stays in that state for the rest of time, we can define the
limiting configuration
$$\tilde{\omega}_{x_1,\dots, x_k}(z)
    =\lim_{t \to \infty} \tilde{\omega}_{x_1,\dots, x_k}(z,t) $$
We say that {\bf mutual unbounded growth} or {\bf coexistence}
occurs if the limiting configuration has infinitely many $z$ in
state i for all $i\leq k$. More precisely we define
$\mug(x_1,\dots,x_k)$ to be the event that
$$|\{z:\tilde{\omega}_{x_1,\dots,x_k}(z)=1\}|
    =\dots    =|\{\tilde{\omega}_{x_1,\dots, x_k}(z)=k\}|
    =\infty.$$
We refer to this event as {\bf coexistence} or {\bf mutual
unbounded growth}.

\subsection{Results}
Our results depend on the geometry of $\ball$. Let $\sides$ be the
number of sides of $\partial \ball$ if $\partial \ball$ is a
polygon and infinity if $\partial \ball$ is not a polygon. Note
that by symmetry we have $\sides \geq 4$ for any good measure
$\mu$. Let $G(x_1,\dots,x_k)$ be the event that there exist
disjoint geodesics $g_i$ starting at $x_i$. In this paper we prove
the following theorem about general first passage percolation.
\begin{theorem} \label{mainfpp}
Let $\fpp$ be good. For any $\epsilon>0$ and $k \leq \sides$ there
exists $x_1, \dots x_k$
 such that
 $$\prob(G(x_1,\dots,x_k))>1-\epsilon.$$
\end{theorem}
We also get the two closely related theorems. Let $A$ finite
subset of $\edges(\Z^2)$.  Let $\{(c_a,d_a)\}_{a \in A}$ be a
collection of intervals with $0<c_a<d_a$ for all $a \in A$. Let
$A'$ be the event that $\omega(a) \in (c_a,d_a)$ for all $a \in
A$. Let $B'$ be any event such that $P(B')>0$ and $B'$ does not
depend on $\omega|_A$ (if $\omega \in B'$ and
$\omega|_{A^c}=\omega'|_{A^c}$ then $\omega' \in B'$). We say that
$\mu$ has {\bf finite energy} if
 $$\P(A' \cap B')>0$$ for all such events $A'$ and $B'$.
\begin{theorem} \label{exponential}
If $\fpp$ is good and has finite energy then for any $k \leq \sides$
$$\prob(|K(\Gamma(\0))| \geq k)>0.$$
\end{theorem}

\begin{theorem} \label{mainfpp2}
Let $\fpp$ be good. For any $k \leq \sides/2$
 $$\prob(|K(\Gamma(\0))|\geq k)=1.$$
\end{theorem}

Theorem \ref{exponential} extends a theorem of \olle and Pemantle
\cite{HP}. They proved that under the same hypothesis that
$$ \prob(|K(\Gamma(\0))|>1)>0.$$
Garet and Marchand \cite{G} and Hoffman \cite{hoffman} extended
the results of \olle and Pemantle to a general class of first
passage percolation processes in any dimension.

As an easy consequence of Theorem \ref{mainfpp2} we get
\begin{cor} \label{followsfromhm}
There exists a good measure $\fpp$ such that
$$\prob(|K(\Gamma(\0))|=\infty)=1.$$
\end{cor}
\begin{proof}
This follows easily from Theorem \ref{mainfpp2} and \cite{H} where
it is proven that there is a good measure $\fpp$ such that $\ball$
is the unit disk.
\end{proof}

Our main result on a multiple color  Richardson's growth model is
that with positive probability coexistence occurs.
\begin{thm} \label{maingrowth} If $\mu$ is good and $k \leq \sides$ then
for any $\epsilon>0$ there exist $x_1,\dots,x_k$ such that
 $$\prob(\mug(x_1,\dots,x_k))>1-\epsilon.$$
\end{thm}
\olle and Pemantle \cite{HP} proved that if $\fpp$ is i.i.d.\ with
exponential distribution then
$$ \prob(\mug((0,0),(0,1)))>0.$$
Garet and Marchand \cite{G} and Hoffman \cite{hoffman} proved that
in any dimension mutual unbounded growth is possible when $k=2$.
Our result extends the previous results in two ways.  First it
shows that coexistence is possible with four colors. It also
shows that the points $x_1,\dots,x_k$ can be chosen such that the probability of coexistence approaches one.  None of
the three proofs that coexistence is possible in the two color
Richardson's growth model were able to show that the probability
of coexistence went to one as the initial sites $x_1$ and $x_2$
moved farther apart.

\begin{corollary}
There exists a nontrivial i.i.d.\ measure $\fpp$ and
$x_1,\dots,x_8$ such that
$$\prob(C(x_1,\dots,x_8))>0.$$
\end{corollary}

\begin{proof}
By \cite{CD} there exists a $\fpp$ which is i.i.d.\ such that
$\partial \ball$ is neither a square nor a diamond. As $\ball$ is
symmetric $\sides \geq 8$. Thus the result follows from Theorem
\ref{maingrowth}.
\end{proof}

\section{Notation}

Much of the notation that we introduce is related to the shape
$\ball$. For $v \in \R^2\setminus (0,0)$ let
$$T^*(v)=\frac{1}{\sup\{k: k v \in \ball \}}. $$
It is not hard to check that $T^*$ is a norm on $\R^2$ whose unit ball is $\ball.$ Thus it is equivalent with Euclidean distance.  It might be helpful to note that Theorem \ref{shape1} implies that
\begin{equation} \label{T}
\lim \frac 1n \tau(\0,nv)=T^*(v) \ \ a.s.
\end{equation}
Also we have that $T^*(\alpha v)=\alpha T^*(v)$ and $T^*(v)=1$ for all $v \in
\partial \ball$ and all $\alpha \in \R$.  We use $||v||=\sqrt{v_1^2+v_2^2}$ to represent
the length of $v$.

Let the set $\good$ consist of all $v \in
\partial \ball$ such that there is a unique line $L_v$ which is
tangent to $\ball$ through $v$. For such a $v$ let $w(v)$ be a
unit vector parallel to $L_v$. Let $L_{n,v}$ be the line through
$nv$ in the direction of $w(v)$. We start with two lemmas about
the set $\good$.

\begin{lemma} \label{marchand}
For any $k \leq \sides$ there exists points
$v_1,\dots,v_k  \in \good$ such that
 the lines $L_{v_i}$ are distinct for all $i$.
\end{lemma}

\begin{proof}
If $\sides<\infty$ then $\partial \ball$ is a polygon and the
lemma is obvious.  For $v$ in the first quadrant define $w^{+}(v)$
to be the largest angle (measured counterclockwise) between the
positive $x$-axis and a line through $v$ that does not intersect
the interior of $\ball$. (At least one such line exists by the
convexity of $\ball$.) Define $w^{-}(v)$ to be the smallest angle
(measured counterclockwise) between the positive $x$-axis and a
line through $v$ that does not intersect the interior of $\ball$.

As $v$ rotates from being parallel to the positive $x$-axis to
being parallel to the positive $y$-axis both $w^+$ and $w^-$ are
nondecreasing. Thus they are continuous almost everywhere.  It is
easy to check that $\partial \ball$ has a unique tangent at $v$ if
and only if $w^+(v)=w^-(v)$. As $\ball$ is convex there is a
unique tangent line at almost every point in $\partial \ball$ and
the two functions are equal for almost every $v$.

If $\partial \ball$ is not a polygon then $w^+(v)$ takes on
infinitely many values.  For every $i \in \N$ choose $v_i$ such
that $w^+(v_{i+1})>w^+(v_i)$ for all $i$.  For each $i$ choose
$v'_i$ such that $v'_i$ is in the arc of $\partial \ball$ from
$v_{2i}$ to $v_{2i+1}$ and $w^+(v'_i)=w^-(v'_i)$. This is possible
because the two functions are equal almost everywhere so there
exists a point of equality on every arc of positive length. Thus
at each $v'_i$ there is a unique tangent line to $\partial \ball$.
For any $i > j$ we have that $$w^+(v'_i)\geq
w^+(v_{2i})>w^+(v_{2j+1})\geq w^+(v'_j)$$ and the tangent lines at
$v'_i$ and $v'_j$ are distinct.
\end{proof}

\begin{lemma} \label{garet}
There is a unique line tangent to $\partial \ball$ at the point
$v$ if and only if
\begin{equation} \label{smooth}
\lim_{b \to 0} \frac{T^*(v+w(v)b)-1}{|b|}=0.
\end{equation}
\end{lemma}

\begin{proof}
Fix $v \in \partial \ball$. For $v'\in \partial \ball$ and $v'$ not
parallel to $w(v)$ we can find $a$ and $b$ such that $v'=av+abw(v)$.
Then we have
$$T^*(v+bw(v))=T^*\left(\frac1a{v'}\right)=\frac1a T^*(v')=\frac1a.$$
As $v'$ approaches $v$ we have $a \to 1$
and $b\to 0.$

It is easy to check that $\partial \ball$ having a unique tangent
line at $v$ is equivalent to
$$\lim_{v' \to v,v' \in \partial \ball}\frac{||v+bw(v)-v' ||}{||bw(v)||} = 0.$$
As
$$\frac{||v+bw(v)-v' ||}{||bw(v)||}=\frac{(\frac1a -1)||v'||}{|b|\cdot||w(v)||}$$
having a unique tangent line is equivalent to
$$\lim_{v' \to v,v' \in \partial \ball}\frac{(\frac1a -1)||v'||}{|b|\cdot||w(v)||} = 0.$$

Since $||v'|| \to ||v||\neq 0$ as $v' \to v$ having a unique
tangent line is equivalent to
$$\lim_{v' \to v,v' \in \partial \ball}\frac{\frac1a -1}{|b|} = 0.$$
We have that $b \to 0$ is equivalent to $v' \to v$ for $v' \in
\partial \ball$. Thus $\partial \ball$ having a unique tangent at $v$ is
equivalent to
$$
 \lim_{v' \to v,v' \in \partial \ball}\frac{\frac1a-1}{|b|} =
 \lim_{b \to 0}\frac{T^*(v+bw(v))-1}{|b|}= 0.$$
\end{proof}

Let $S \subset \R^2$ we define the function
$$B_{S}(x,y)=\inf_{z \in S}\tau(x,z)-\inf_{z \in S}\tau(y,z).$$

\begin{lemma} \label{busemannproperties}
For any set $S\subset \z^2$ and any $x,y,z \in \z^2$
\begin{enumerate}
\item $B_S(x,y) \leq \tau(x,y)$ and
\item $B_S(x,y)+B_S(y,z)=B_S(x,z)$.
\end{enumerate}
\end{lemma}

\begin{proof}
These properties follow easily from the subadditivity of $\tau$
and the definition of $B_S$.
\end{proof}

These functions are useful in analyzing the growth model because
of the following fact.

\begin{lemma} \label{mugconditions}
If there exists $c>0$, and $x_1,\dots,x_k \in V$ such that
$$\prob(B_{L_{n,v_i}}(x_j,x_i)>0 \ \forall i \neq j)\geq 1-c$$
for infinitely many $n$ then
$$\prob(C(x_1,\dots,x_k))\geq 1-c.$$
\end{lemma}

\begin{proof}
If for a fixed $i$ and all $j \neq i$
$$B_{L_{n,v_i}}(x_j,x_i)>0$$
 then there exists $z \in L_{n,v_i}$ such that
$\tau(z,x_i)<\tau(z,x_j)$ for all $j\neq i$. Thus there is a $z\in
L_{n,v_i}$ such that $$\tilde \omega_{x_1,\dots,x_k}(z)=i.$$ For a
fixed $i$ each $z$ is in only one $L_{n,v_i}$, so if there exist
infinitely many $n$ such that for all $i$ and $j \neq i$
$$B_{L_{n,v_i}}(x_j,x_i)>0$$
then for every $i$ there are infinitely many $z$ such that
$$\tilde \omega_{x_1,\dots,x_k}(z)=i.$$

By assumption we have that there exist infinitely many $n$ such
that
$$\prob(B_{L_{n,v_i}}(x_j,x_i)>0 \text{ for all $i$ and $j \neq i$})\geq 1-c.$$
Thus  we have that with probability at
least $1-c$ there exist infinitely many $n$ such that for all $i$
and $j \neq i$
$$B_{L_{n,v_i}}(x_j,x_i)>0.$$
In conjunction with  the previous paragraph this proves the lemma.
\end{proof}


\begin{lemma} \label{mugimpliesgeodesics}
$$\prob(G(x_1,\dots,x_k)) \geq \prob(C(x_1,\dots,x_k)).$$
\end{lemma}

\begin{proof}
For any
$$y \in \{z: \tilde\omega_{x_1,\dots,x_k}(z)=i\}$$
the geodesic from $x_i$ to $y$ lies entirely in
$$\{z: \tilde\omega_{x_1,\dots,x_k}(z)=i\}.$$
 By compactness if
 $$|\{z: \tilde\omega_{x_1,\dots,x_k}(z)=i\}|=\infty$$
then there exists an infinite geodesic $g_i$ which is contained in
the vertices
$$\{z: \tilde\omega_{x_1,\dots,x_k}(z)=i\}.$$
\end{proof}

\begin{lemma} \label{reallysmall}
For all $v \in \R^2\setminus \0$ and $\epsilon>0$ there exists
$\delta=\delta(\epsilon,v)>0$ and $M_0=M_0(\epsilon,v)$ such that
for all  $M>M_0$, all events $E$ with $\prob(E)<\delta$ and any
$r\in\R^2$ we have that
$$\Exp(\tau(r,r+Mv)\indicator_E)<M\epsilon.$$
\end{lemma}

\begin{proof}
Consider the space $\bar {\Omega}=[0,\infty)^{\edges(\zd)} \times
[0,1) \times [0,1)$ with measure $\bar \mu$ the direct product of
$\mu$ with Lebesgue measure.  We write $\bar \prob$ and
$\Exp_{\bar \mu}$ for probability and expectation with respect to
$\bar \mu$. Any vector $v \in \R^2 \setminus 0$ acts on $\bar
\Omega$ in the following manner.

For any $v \in \R^2 \setminus \0$ and $(\omega,a,b) \in \bar
\Omega$ we have
$$\bar{\sigma}_v(\omega,a,b)=(\sigma_{v'}(\omega),c,d)$$ where
$v+(a,b)=v'+(c,d)$, $v' \in \z^2$ and $c,d \in [0,1) \times
[0,1)$.  For convenience we often write for $(a,b) \in \R^2$
$$(\omega,a,b)=(\sigma_{v'}(\omega),c,d)$$ where $(a,b)=v'+(c,d)$,
$v' \in \z^2$ and $c,d \in [0,1) \times [0,1)$. For any $a,b,c,d
\in \R$ and $\omega$ we write
 $$\tau((\omega,(a,b)),(\omega,(c,d)))=\tau^{\omega}((a,b),(c,d)).$$
We also define the function
$$f(\omega,a,b)=\tau((\omega,(a,b)),(\omega,(a,b)+v)).$$ Note that $f$ is in
$L^1$.

For any set $E $ with $\prob(E)< \delta$ we write $\bar E=E \times
[0,1) \times [0,1) $ and we have $\bar \prob(E)<\delta$. For any
$M \in \R$ choose $k$ such that $k\leq M \leq k+1$.

\begin{eqnarray*}
\lefteqn{\frac1M \Exp(\tau(r,r+Mv))\indicator_{E}}\\
 & \leq& \frac1k \left(\Exp \sup_{a,b\in[0,1)} \tau ((a,b),kv+(a,b))\indicator_E
     +\Exp\sup_{M \in [k,k+1)}(\tau(r+kv,r+Mv))\right)\\
 & \leq& \frac1k \bigg(\Exp_{\bar \mu} \sup_{a,b\in[0,1)} \tau ((\omega,(a,b)),(\omega,kv+(a,b)))\indicator_{\bar E}
    +\Exp\sup_{M \in [k,k+1)}(\tau(r+kv,r+Mv)) \bigg)\\
 & \leq& \frac1k \bigg(\Exp_{\bar \mu}  \tau ((\omega,(a,b)),(\omega,kv+(a,b)))\indicator_{\bar E}
    +\Exp \sup_{a,b,c,d\in[0,1)}\tau((a,b),(c,d))\\
&&  +\Exp \sup_{a,b,c,d\in[0,1)}\tau(kv+(a,b),kv+(c,d))
    +\Exp \sup_{M \in [k,k+1)}(\tau(r+kv,r+Mv))
\bigg)\\
 & \leq& \frac1k \bigg(\Exp_{\bar \mu}  \tau ((\omega,(a,b)),(\omega,kv+(a,b))\indicator_{\bar E}
    +\Exp \sup_{a,b,c,d\in [0,1)}\tau((a,b),(c,d))\\
  &&  +\Exp \sup_{a,b\in[0,1)}\tau(kv,kv+(a,b))
+\Exp\sup_{M \in [k,k+1)}(\tau(r+kv,r+Mv)) \bigg).
\end{eqnarray*}
The second, third and fourth terms in the last inequality are
bounded independent of $k$. Thus their contribution to the right
hand side goes to zero as $k$ goes to infinity.
Then we have that
\begin{eqnarray}
\frac1k \tau((\omega,(a,b)),(\omega,kv+(a,b))) &\leq& \frac1k
\sum_{0}^{k-1} \tau (\omega,jv+(a,b),\omega,(j+1)v+(a,b))\nonumber\\
    &\leq& \frac1k \sum_{0}^{k-1} f(\bar{\sigma}_v^j(\omega,a,b)).\label{suspension}
    \end{eqnarray}
By the ergodic theorem the sum on the righthand side of
(\ref{suspension}) is converging to an $L^1$ function almost
everywhere and in $L^1$.   Thus we can choose $\delta$ such that
$\bar{\prob}(\bar E)<\delta$ implies
$$\frac1k \Exp_{\bar \mu}  \tau ((\omega,(a,b)),(\omega,kv+(a,b)))\indicator_{\bar E}
    \leq \Exp_{\bar \mu} \left( \frac1k \sum_{0}^{k-1} f(\bar{\sigma}_v^j(\omega,a,b)) \indicator_{\bar E}\right)
  <\epsilon.$$
 This proves the lemma.
\end{proof}


We use this lemma in two contexts.

\begin{cor} \label{small}
For all $v \in \R^2\setminus \0$, $r\in \R^2$, $\epsilon>0$ and $M
\in\R$ let $E=E(M,v,r,\epsilon)$ be the event that
$$\tau(r,r+Mv)>(1+\epsilon/2)M.$$ There exists
$M_0=M_0(v,r,\epsilon)$ such that for all $M>M_0$
$$\E\left( \tau(r,r+Mv)\indicator_E\right)\leq M\epsilon.$$
\end{cor}

\begin{proof}
Fix $v$, $r$, $m$ and $\epsilon$.
 By Theorem \ref{shape1} we have
that $\P(E)\to 0$ as $M \to \infty$.  Thus we can apply Lemma
\ref{reallysmall} 
 to prove the corollary.

\end{proof}

\begin{cor} \label{smaller}
For any $v \in V$, $r \in \R^2$, $m \in \R$ and $\epsilon>0$ there
exists $M_0=M_0(v,r,m,\epsilon)$ such that for all $M>M_0$
$$\E\left( \tau(r-mv,r+Mv)\right)\leq M(1+\epsilon).$$
\end{cor}

\begin{proof}
Fix $v$, $r$, $m$ and $\epsilon$.   Let $E=E(M)$ be the event that
$\tau(r-mv,r+Mv)>M(1+\epsilon/2).$  By Theorem \ref{shape1} we
have that $\P(E)\to 0$ as $M \to \infty$.  Thus we can apply Lemma
\ref{reallysmall}
  to prove
$$\E(\tau(r-mv,r+Mv)\indicator_{E})\leq M(\epsilon/2).$$
By the definition of $E$
$$\E(\tau(r-mv,r+Mv)\indicator_{E^c})\leq M(1+\epsilon/2).$$
Putting those two together proves the corollary.
\end{proof}

%



\section{Outline}

We start by outlining a possible method to prove that there are
infinitely many geodesics starting at the origin.  Then we show
the portion of this plan that we can not prove.  Finally we show
how to adapt this method to get the results in this paper.

It is easy to construct geodesics beginning at $\0$.  We can take
any sequence $W_1,W_2,\dots$ of disjoint subsets of $\z^2$ and
consider $G(\0,W_n)$, the geodesic from $\0$ to $W_n$.  (The finite geodesic $G(0,W_n)$ is well defined a.s.\ because the measure $\mu$ is good so $R$ is bounded. Then Theorem \ref{shape1} implies the existence of the finite geodesic.) Using
compactness it is easy to show that there exists a subsequence
$n_k$ such that $G(\0,W_{n_k})$ converges to an infinite geodesic.

If we take two sequences of sets $W_n$ and $W'_n$ we can construct
a geodesic for each sequence.  It is difficult to determine
whether or not the two sequences produce the same or different
geodesics. The tool that we use to distinguish the geodesics are
Busemann functions.  Every geodesic generates a Busemann function
as follows.

For any $x,y \in \z^2$ and infinite geodesic
$G=(v_0,v_1,v_2,\dots)$ we can define
$$\hat B^{\omega}_{G}(x,y)= \hat B_G(x,y)=\lim_{n \to \infty} \tau(x,v_n)-\tau(y,v_n).$$
To see the limit exists first note that
\begin{eqnarray*}
\hat B_G(x,y)& =& \lim_{n \to \infty} \tau(x,v_n)-\tau(y,v_n)\\
        & =& \lim_{n \to \infty} \tau(x,v_n)-\tau(v_0,v_n)+\tau(v_0,v_n)-\tau(y,v_n)\\
        & =& \lim_{n \to \infty} (\tau(x,v_n)-\tau(v_0,v_n))+\lim_{n \to \infty}(\tau(v_0,v_n)-\tau(y,v_n)).
\end{eqnarray*}
As $G$ is a geodesic the two sequences in the right hand side of
the last line are bounded and monotonic so they converge. Thus
$\hat B_G(x,y)$ is well defined.

Two distinct geodesics may generate the same Busemann function but
distinct Busemann functions mean that there exist distinct
geodesics.

To construct a geodesic we pick $(a,b) \in \z^2$ and we set $W_n$
to be
$$W_n=\{w\in \z^2:w\cdot (a,b) \geq n\}.$$
If we could show that for every $z\in \z^2$ that $G_n(z)$, the
geodesic from $z$ to $W_n$, converges then it would be possible to show
that
\be \label{star1}
\lim_{M \to \infty}\frac{1}{M}B(\0,(bM,-aM))=0.
\ee
and
\be \label{star3}
\lim_{M \to \infty}\frac{1}{M}B(\0,(a M,-b M))
    =\inf_{v\cdot(a,b)=a^2+b^2}T^*(v).
\ee
Thus for any $(a',b')$ which is not a scalar multiple of $(a,b)$
we would be able to show that
\be \label{star2}
\lim_{M \to \infty}\frac{1}{M}B(\0,(b'M,-a'M))\neq 0.
\ee
Thus for any $(a,b)$ and $(a',b')$ which are not scalar multiples
we get distinct geodesics.  In this way it would be possible to
construct an infinite sequence of distinct geodesics. We are
unable to show that the geodesics $G(z,W_n)$ converge. But for
some $(a,b) \in \z^2$ we can establish  versions of (\ref{star1})
and (\ref{star3}).  These are Lemmas \ref{slope} and \ref{B_n}.
These lemmas form the heart of our proof.

\section{Proofs}

Although it is convenient to write $\tau(x,y)$ for $x,y \in \R^2$,
the distribution of $\tau(x,y)$ is equal to the distribution of
$\tau(x+z,y+z)$ only if $z \in \z^2$.  For $z \not\in \z^2$ the
distribution of $\tau(x,y)$ may not be equal to the distribution
of $\tau(x+z,y+z)$ which will make the notation more complicated.
But the distributions are close enough so that this lack of shift
invariance for noninteger translations will not cause any
significant problems.  To deal we this lack of translation
invariance we let
$$I(a,b)=I_v(a,b)=\sup_{x \in L_{a,v}}\left(\Exp \left( \tau(x,L_{b,v})\right)\right)$$
In the next three lemmas we show
$$\sup_{x \in L_{a,v}}\left|I(a,b)-\Exp \left(
\tau(x,L_{b,v})\right)\right|$$ is bounded uniformly in $a,b$ and
$v$.
\begin{lemma} \label{anal}
For any $v \in V$ and $u \in \R^2\setminus 0$ such that
 $u =\alpha v + \gamma w(v)$
 $$\frac{\alpha||v||}{||u||} <\sqrt{2}.$$
\end{lemma}

\begin{proof}
First we show that for any $v \in V$ that lies in the first octant
(between the lines $y=0$ and $y=x$ with $x>0$) that $w(v)$ then
points in one of the octants between the lines $x=0$ and $x=-y$.

Let $\tilde v$ be the image of $v$ under reflection about the line
$x=y$ and $v^*$ be the image of $v$ under reflection about the
line $x=0$.  The line from $v$ to $\tilde v$ is parallel to the
line $x=-y$ while the line from $v$ to $v^*$ is parallel to the
line $x=0$.


By the convexity of $\ball$ we have that for any two points in
$\ball$ and any line tangent to $\partial \ball$ the two points
lie on the same side of the line (or in one closed halfplane).
Thus $\0,$ $v^*$ and $\tilde v$ all lie on the same side of
$L_{1,v}$. This implies $L_{1,v}$ does not intersect the interior
of the line segment between $v^*$ and $\tilde v$ and $w(v)$ points
in the octants between the lines $x=0$ and $x=-y$.

Then
\begin{equation} \label{dot} \nonumber
|v\cdot w(v)|\leq \frac{\sqrt{2}}{2}||v||\cdot||w(v)||
\end{equation}
and the angle between $v$ and $w(v)$ is at least 45 degrees.  By
the symmetry of $\ball$ this inequality holds for all $v\in V$.
For a fixed $\alpha$ and $v$ the value of $\gamma$ which minimizes
$||u||$ occurs when the points $\0$, $\alpha v$ and $u$ form a
right triangle.  As the angle between $v$ and $w(v)$ is at least
45 degrees we have that
$$\frac{\alpha ||v||}{||u||}<\sqrt{2}.$$

\end{proof}

\begin{lemma} \label{white}
Let $v \in V$, $x_1,x_3 \in \R^2$ and $n_1,n_2,n_3,n_4 \in \R$
with $x_1 \in L_{n_1,v}$, $x_3 \in L_{n_3,v}$, $n_1<n_2$,
$n_3<n_4$ and
$$(n_2-n_1)-(n_4-n_3)\geq 2/||v||.$$
Then
$$\Exp(\tau(x_1,L_{n_2,v}))\geq \Exp(\tau(x_3,L_{n_4,v})).$$
 This implies that for any $m>2/||v||$ and any $r \in L_{n_1,v}$
\begin{equation}\label{weaver}
\Exp(\tau(r-mv,L_{n_2,v})) \geq I(n_1,n_2).
\end{equation}
\end{lemma}

\begin{proof}
First we define $\tilde x_1$ and $\tilde x_3$ to be the points in
$\Z^2$ closest to $x_1$ and $x_3$ respectively (i.e. $x_1 \in
\tilde x_1 +[-1/2,1/2)^2$ and $x_3 \in \tilde x_3 +[-1/2,1/2)^2$).

Next define $\alpha_1,\alpha_3,\gamma_1,\gamma_3 \in \R$ such that
$$\tilde x_1-x_1= \alpha_1 v +\gamma_1 w(v) \text{ and }
\tilde x_3-x_3= \alpha_3 v +\gamma_3 w(v).$$ By the definition of
the $\tilde x_i$ we have that $||x_i-\tilde x_i||\leq \sqrt{2}/2$.
Thus by Lemma \ref{anal} we have that $|\alpha_1|,|\alpha_3| \leq
1/||v||.$

Then define $\tilde \alpha$ and $\tilde \gamma$ such that
\begin{eqnarray*}
\tilde \alpha v+\tilde \gamma w(v)
 &=& \tilde x_3 -\tilde x_1\\
 &=& x_3 - x_1 +(\alpha_3-\alpha_1)v+(\gamma_3-\gamma_1)w(v)\\
 &=& (n_3-n_1 +\alpha_3-\alpha_1)v+Cw(v)
\end{eqnarray*}
for some $C \in \R$.  Also
$$|\tilde \alpha -(n_3-n_1)| =|\alpha_3-\alpha_1| \leq 2/||v||$$
or
$$(n_3-n_1)-2/||v|| \leq \tilde \alpha. $$
Then
$$L_{n_2,v}+(\tilde x_3 -\tilde x_1)=L_{n_2+\tilde \alpha,v}.$$
Also
$$\Exp(\tau(x_1,L_{n_2,v}))=\Exp(\tau(\tilde x_1,L_{n_2,v}))
 =\Exp(\tau(\tilde x_3,L_{n_2+\tilde \alpha,v}))
 =\Exp(\tau(x_3,L_{n_2+\tilde \alpha,v})).$$
 The first and third inequalities are due to the definition of
 $\tilde x_1$ and $\tilde x_3$ respectively while the second is
 due to the shift invariance of the distribution under shifts in
 $\Z^2$.  (The image of $L_{n_2,v}$ under translation by $\tilde x_3 -\tilde x_1$ is $L_{n_2+\tilde \alpha,v}$.) Thus
\begin{equation}\label{christmas}
 \Exp(\tau(x_1,L_{n_2,v})) =\Exp(\tau(x_3,L_{n_2+\tilde \alpha,v}))\geq \Exp(\tau(x_3,L_{n_4,v}))
\end{equation}
if and only if $n_4\leq n_2+\tilde \alpha$. As
\begin{eqnarray*}
2/||v|| &\leq&(n_2-n_1)-(n_4-n_3) \\
0       &\leq&(n_2-n_4)+(n_3-n_1)-2/||v||\\
0       &\leq&(n_2-n_4)+\tilde \alpha\\
n_4     &\leq&n_2+\tilde \alpha.
\end{eqnarray*}
Thus by (\ref{christmas})
\begin{equation} \label{easter}
 \Exp(\tau(x_1,L_{n_2,v})) \geq \Exp(\tau(x_3,L_{n_4,v}))
\end{equation}
  and the first part of the lemma is true.

For the second statement for any $n_3 \in \R$ take any $r$ and $w$
in $L_{n_1,v}$.
 Apply (\ref{easter}) with  $x_1=w$, $x_3=r-mv$,$n_3=n_1-\alpha$ and $n_2=n_4$ to get
$$ \Exp(\tau(r-mv,L_{n_2,v}))=\Exp(\tau(r-mv,L_{n_4,v})) \geq \Exp(\tau(w,L_{n_2,v})).$$
As this holds for all $w \in L_{n_1,v}$ we have
$$ \Exp(\tau(r-mv,L_{n_2,v})) \geq \sup_{w \in L_{n_1,v}} \Exp(\tau(w,L_{n_2,v}))=I(n_1,n_2).$$
\end{proof}

\begin{lemma} \label{jellystone}
There exists $\beta \in \R$ such that for all $v \in V$ and
$n_1,n_2 \in \R$ with $n_2-n_1>2/||v||$ and for all $w,y \in
L_{n_1,v}$
$$|\Exp(\tau(w,L_{n_2,v}))-\Exp(\tau(y,L_{n_2,v}))|<\beta.$$
 We also have that for any $n_1,n_2,\alpha \in \R$
\begin{equation}\label{newshift}
|I(n_1,n_2)-I(n_1+\alpha,n_2+\alpha)|<\beta.
 \end{equation}
\end{lemma}

\begin{proof}
Pick $\beta$ such that for all $x,z$ with $||x-z||=2$ we have
$\Exp(\tau(x,z)) <\beta.$ Define $x$ and $z$ by $x=y-2v/||v||$ and
$z=y+2v/||v||.$  By Lemma \ref{white} we have
$$\Exp(\tau(z,L_{n_2,v}))\leq \Exp(\tau(y,L_{n_2,v})),\Exp(\tau(w,L_{n_2,v})) \leq \Exp(\tau(x,L_{n_2,v})). $$
Thus
$$0\leq |\Exp(\tau(y,L_{n_2,v}))-\Exp(\tau(w,L_{n_2,v}))|
 \leq \Exp(\tau(x,L_{n_2,v}))-\Exp(\tau(z,L_{n_2,v}))
 \leq \Exp(\tau(x,z)) < \beta.$$

For the second part choose $z \in L_{n_1+\alpha,v}$ and
 $r \in L_{n_1,v}.$  Also choose $m>2/||v||$ such that
 $\E(\tau(r-mv,r))<\beta.$
 There exists $\tilde r \in \Z^2$ and $r-mv \in \tilde r +[1/2,1/2)^2$. Let
$\tilde n_1$ be such that $\tilde r \in L_{\tilde n_1,v}.$
 There exists $\tilde z \in \Z^2$ and $z \in \tilde z +[1/2,1/2)^2$. Let
$\hat n_1$ be such that $\tilde r \in L_{\hat n_1,v}.$
 Then we have $n_2-\tilde n_1>n_2+\alpha - \hat n_1.$
 This implies
 $$\E(\tau(r-mv,L_{n_2,v}))=\E(\tau(\tilde r,L_{n_2,v}))
   >\E(\tau(\tilde z,L_{n_2+\alpha,v}))=\E(\tau(z,L_{n_2+\alpha,v})).$$
As this holds for all $z \in L_{n_1+\alpha,v}$ we have
 \begin{equation}\label{greg}
 \E(\tau(r,L_{n_2,v}))\geq I(n_1+\alpha,n_2+\alpha).
 \end{equation}
As $\E(\tau(r-mv,r))<\beta$ we also have
\begin{equation}\label{oden}
 I(n_1,n_2)+\beta \geq
 \E(\tau(r,L_{n_2,v}))+\beta>\E(\tau(r-mv,L_{n_2,v})).
 \end{equation}
Thus combining (\ref{greg}) and (\ref{oden})
$$I(n_1,n_2)+\beta >\E(\tau(r-mv,L_{n_2,v}))\geq
  I(n_1+\alpha,n_2+\alpha).$$
An analogous argument gives
$$I(n_1+\alpha,n_2+\alpha)+\beta >  I(n_1,n_2)$$
which completes the proof.
\end{proof}

Now we show that for a typical choice of $n,M\in \R$, $v\in V$ and
$r \in \R^2$ we have that $B_{L_{n,v}}(r,r+Mv)$ is close to
$\tau(r,r+Mv)$ (which is close to $M$  because $v \in \partial
\ball$). This (along with Lemma \ref{slope}) is one of two key
steps in showing that for distinct $v,v' \in V$ we will get
distinct \booze functions.

 We define the lower density of
$A \subset \N$ to be
$$\underline{\dens(A)}=\liminf_{N \to \infty} \frac1N |A\cap [1,2,\dots,N]|.$$
Similarly we define
$$\overline{\dens(A)}=\limsup_{N \to \infty} \frac1N |A\cap [1,2,\dots,N]|.$$
We will assume the reader is familiar with all of the normal properties of density of sets, e.g
$$\underline{\dens(A)}+\overline{\dens(A^c)}=1$$
and
 $$\overline{\dens(A \cup B)}\leq \overline{\dens(A)}+\overline{\dens(B)}.$$
We often shorten lower density to density as it will not cause confusion.

\begin{lemma} \label{B_n}
For any $v \in \good$, any $\epsilon>0$, there exists $M_0=M_0(\epsilon,v)$ such that for all $M>M_0$ and all $r\in \R^2$ the density of $n$ such that
\begin{equation} \label{star1''}
 \prob\bigg(M(1-\epsilon) < B_{L_{n,v}}(r,r+Mv)< M(1 + \epsilon)\bigg)
   >1-\epsilon.
\end{equation}
is at least $1-\epsilon$.
\end{lemma}
\begin{proof}
By Lemma \ref{busemannproperties} for any $r,n,M$ and $v$
\be \label{upper1}
B_{L_{n,v}}(r,r+Mv)\leq \tau(r,r+Mv)
\ee
and by Theorem \ref{shape1} for any $r, v$ and sufficiently large
$M$
\be \label{upper}
\prob\bigg(\tau(r,r+Mv)<M(1+\epsilon)\bigg)>1-\epsilon.
\ee
Thus for sufficiently large $M$ the upper bound on
$B_{L_{n,v}}(r,r+Mv)$ is satisfied for all $n$ with probability at
least $1-\epsilon.$

Now we bound the probability that $B_{L_{n,v}}(r,r+Mv)$ is too
small. Let $d$ be such that $r \in L_{d,v}$ and let $m \in \R$ be such that $m||v||>2$. For any sufficiently
large $M$ and any $n\geq d+M$
\begin{eqnarray}
 \lefteqn{I(d,n)-I(d+M,n)} \hspace{1in}&& \nonumber \\
    &\leq& \Exp \left(\inf_{y \in L_{n,v}} \tau(r-mv,y)\right)
        -\sup_{x \in L_{d+M,v}}\left(\Exp \left(\inf_{y \in L_{n,v}} \tau(x,y)\right)\right) \label{I1}\\
    &\leq& \Exp \left(\inf_{y \in L_{n,v}} \tau(r-mv,y)\right)
        -\Exp \left(\inf_{y \in L_{n,v}} \tau(r+Mv,y)\right) \label{I2}\nonumber\\
    &\leq& \Exp \left(\inf_{y \in L_{n,v}} \tau(r-mv,y)
        -\inf_{y \in L_{n,v}} \tau(r+Mv,y)\right) \nonumber\\
    &\leq& \Exp\bigg(B_{L_{n,v}}(r-mv,r+Mv)\bigg) \label{I4} \\
    &\leq& \Exp\bigg(\tau(r-mv,r+Mv)\bigg)  \label{why}\\
    &\leq&  M(1+\epsilon). \label{why2}
\end{eqnarray}
(\ref{I1}) follows from (\ref{weaver}) in Lemma \ref{white} and
the definition of $I(d+M,n)$, (\ref{I4}) follows from the
definition of $B_{L_{n,v}}$, (\ref{why}) follows from Lemma
\ref{busemannproperties}, and (\ref{why2}) follows from Corollary
\ref{smaller}.

  Let $k$ be such that
$$d+kM\leq n <d+(k+1)M.$$
For $k$ and $M$ by Theorem \ref{shape1} and (\ref{newshift})
\begin{eqnarray}
(k+1)M(1-\epsilon)
 &\leq &I(d,n)\\
(k+1)M(1-\epsilon)
 &\leq & I(d,n) +\left(\sum_{l=1}^{k} -I(d+lM,n)+I(d+lM,n)\right)
  \nonumber\\
(k+1)M(1-\epsilon)
 &\leq & \left(\sum_{l=0}^{k-1} I(d+lM,n)-I(d+(l+1)M,n)\right)
 +I(d+kM,n) \nonumber\\
(k+1)M(1-\epsilon)
 &\leq & \left(\sum_{l=0}^{k-1} I(d,n-lM)-I(d+M,n-lM)+2\beta\right)
 +I(d,n-kM). \hspace{.25in}  \label{prince} \nonumber \\
(k+1)M(1-2\epsilon)
 &\leq & \left(\sum_{l=0}^{k-1} I(d,n-lM)-I(d+M,n-lM)\right)
 +I(d,n-kM). \hspace{.25in} \label{newsum} \nonumber \\
(k+1)M(1-2\epsilon)
 &\leq & \left(\sum_{l=0}^{k-1} I(d,n-lM)-I(d+M,n-lM)\right)
 +I(d,d+M). \hspace{.25in} \label{newsum2}
\end{eqnarray}
By  (\ref{why2}) the sum  in the right hand side of (\ref{newsum2})
is the sum of $k+1$ terms bounded above by $M(1+\epsilon)$.
Thus the number of
$l<k$ such that
 $$ I(d,n-lM)-I(d+M,n-lM)> M(1-\sqrt{\epsilon})$$
is at least $k(1-4\sqrt{\epsilon}).$ The above result held for all $\epsilon>0$ and all $M=M(\epsilon)$ sufficiently large. Thus we get that for any $\epsilon>0$ and any $M \in \N$
sufficiently large and any $j \in [0,1,2,\dots, M-1]$ the density of $n$ such that
\be \label{fuji} I(d,j+Mn)-I(d+M,j+Mn)> M(1-\epsilon)\ee
is at least $1-\epsilon$.  Combining this result for all $j\in [0,1,2,\dots, M-1]$
we get that for any $\epsilon>0$ and any $M \in \N$
sufficiently large the density of $n$ such that
\be \label{fuji2} I(d,n)-I(d+M,n)> M(1-\epsilon)\ee
is at least $1-\epsilon$.

Now we show that for any $M,n,r$ and $v$ such that (\ref{fuji2}) is
satisfied we have that with high probability $B_{L_{n,v}}(r,r+Mv)$
is large. Let $E$ be the event that
$$ B_{L_{n,v}}(r,r+Mv)>M(1+\epsilon/2).$$
By Lemma \ref{busemannproperties} and Theorem \ref{shape1} we can
make $\prob(E)$ arbitrarily small by making $M$ sufficiently
large.  Then we get
\begin{eqnarray}
\Exp\bigg(B_{L_{n,v}}(r,r+Mv) \bigg)   &=& \Exp\bigg(\inf_{z\in L_{n,v}}\tau(r,z)\bigg)-
\Exp\bigg(\inf_{z\in L_{n,v}}\tau(r+Mv,z)\bigg)\nonumber\\
\Exp\bigg(B_{L_{n,v}}(r,r+Mv) \bigg)   &>& I(d,n)-\beta-I(d+M,n) \hspace{.25in}\label{peyton}\\
\Exp\bigg(B_{L_{n,v}}(r,r+Mv)\bigg) &>&  M(1-\epsilon)-\beta \label{ochocinco} \\
\Exp\bigg(B_{L_{n,v}}(r,r+Mv)\bigg)&>&  M(1-2\epsilon)
  \label{bad30}\\
\Exp\bigg(B_{L_{n,v}}(r,r+Mv)\indicator_{E}\bigg)&&\nonumber\\
+\Exp\bigg(B_{L_{n,v}}(r,r+Mv)\indicator_{E^C}\bigg)&>&
M(1-2\epsilon)
  \label{bad2}\\
 \Exp\bigg(B_{L_{n,v}}(r,r+Mv)\indicator_{E^C}\bigg)
 &>& M(1-2\epsilon) - \Exp\bigg(B_{L_{n,v}}(r,r+Mv)\indicator_{E}\bigg)   \nonumber \\
 \Exp\bigg(B_{L_{n,v}}(r,r+Mv)\indicator_{E^C}\bigg)
 &>& M(1-2\epsilon) - \Exp(\tau(r,r+Mv)\indicator_{E}) \label{bad3}   \\
 \Exp\bigg(B_{L_{n,v}}(r,r+Mv){\indicator}_{E^C}\bigg)
 &>& M(1-3\epsilon). \label{bad4}
\end{eqnarray}
(\ref{peyton}) follows from Lemma \ref{jellystone}.
(\ref{ochocinco}) follows from (\ref{fuji}). (\ref{bad30}) holds
for large $M$.
(\ref{bad3}) follows from Lemma \ref{busemannproperties} and
Theorem \ref{shape1}. (\ref{bad4}) is due to Corollary
\ref{small}.

As the expected value of the function
$$B_{L_{n,v}}(r,r+Mv)\indicator_{E^C}$$
is close to its maximum, $M(1+\epsilon),$ we get that with high
probability the function is close to its maximum. Thus we get that
for any $\epsilon>0$ (possibly larger than the previous $\epsilon$ but still arbitrarily small) and all sufficiently large $M$, the set of
$n$ such that
\be \label{lower}
\prob\bigg(M(1-\epsilon)<B_{L_{n,v}}(r,r+Mv)\bigg)
   >1-\epsilon
   \ee
has density at least $1-\epsilon$. Putting together
(\ref{upper1}), (\ref{upper}) and (\ref{lower}) proves the lemma.
\end{proof}

\begin{lemma} \label{slope}
For any $v \in \good$,  $\epsilon>0$, there exists $M_0=M_0(\epsilon,v)$ such that for any $M\in \R$ with
$|M|>M_0$ and any $r\in \R^2$ the density
of $n$ such that
$$\prob\bigg(\left|B_{L_{n,v}}(r,r+M w(v))\right| < \epsilon |M|\bigg)>1-\epsilon$$
is at least $1-\epsilon$.
\end{lemma}

\begin{proof}
First we prove the upper bound in the case that $M$ is positive.
Fix $\epsilon>0$. Since $\ball$ has a unique tangent line at $v$
by (\ref{smooth}) we can find $b>0$ such that
\begin{equation} \label{star0'}
T^*(v+bw(v))(1+\epsilon b)<(1+2\epsilon b).
\end{equation}
By Lemma \ref{B_n} for any $\epsilon,b>0$ and all sufficiently
large $M$ the density of $n$ such that
\begin{equation} \label{star1'}
\prob\bigg(B_{L_{n,v}}(r-Mv,r)\geq M(1 - \epsilon b)\bigg)
   >1-\epsilon
\end{equation}
is at least $1-\epsilon$. By Theorem \ref{shape1} for any
$\epsilon,b>0$ and all sufficiently large $M$
\begin{equation} \label{star2'}
\prob\bigg(\tau(r-Mv,r+Mbw(v))\leq MT^*(v+bw(v))(1 + \epsilon b) \
   \bigg)
   >1-\epsilon.
\end{equation}
Choose $M$ large enough such that both (\ref{star1'}) and
(\ref{star2'}) are satisfied.

Thus with probability at least $1-2 \epsilon$ the density of $n$
such that the following inequalities are satisfied is at least
$1-\epsilon.$
\begin{eqnarray}
\tau(r-Mv,r+Mbw(v))         & \geq & B_{L_{n,v}}(r-Mv,r+Mbw(v)) \label{q1} \nonumber\\
\tau(r-Mv,r+Mbw(v))         & \geq & B_{L_{n,v}}(r-Mv,r) +  B_{L_{n,v}}(r,r+bMw(v))\label{q2} \nonumber \\
\tau(r-Mv,r+Mbw(v)) - B_{L_{n,v}}(r-Mv,r)
                         & \geq &  B_{L_{n,v}}(r,r+bMw(v))\label{q3} \nonumber \\
M T^*(v+bw(v))(1+\epsilon b)- M(1-\epsilon b)
                         & \geq &  B_{L_{n,v}}(r,r+bMw(v)) \label{q4}\\
M (1+2\epsilon b) - M(1-\epsilon b)
                         & \geq &  B_{L_{n,v}}(r,r+bMw(v))) \label{q5}\\
3bM \epsilon
                         & \geq &  B_{L_{n,v}}(r,r+bMw(v))).
                         \label{q6} \nonumber
\end{eqnarray}
The first two lines follow deterministically from Lemma
\ref{busemannproperties}.  (\ref{q4}) is true with probability at
least $1-2\epsilon$.  This follows from (\ref{star2'})  and
(\ref{star1'}).  (\ref{q5}) follows from (\ref{star0'}). Thus we
have that for any sufficiently large $M$ the density of $n$ such
that
$$\prob\bigg( B_{L_{n,v}}(r,r+bMw(v))  \leq  3bM\epsilon \bigg)
   >1-2\epsilon$$
is at least $1-\epsilon$.  By replacing $w(v)$ with $-w(v)$ and
interchanging $r$ and $r+bMw(v)$ we get that for any sufficiently
large $M$ the density of $n$ such that
$$\prob\bigg( B_{L_{n,v}}(r,r+bMw(v))  \geq - 3bM\epsilon \bigg)
   >1-2\epsilon$$
is at least $1-\epsilon$. The case that $M$ is negative follows in
the same manner by replacing $w(v)$ with $-w(v)$. As $\epsilon$
was arbitrary the lemma follows.
 \end{proof}


\begin{lemma} \label{lem:highprob}
Let $v \in V$. For all $y \in \R^2$ let $s=s(v,y)$ and $t=t(v,y)$
be such that
\be v+sw(v)=y+tv.\label{regan}\ee Then for all
$\epsilon>0$ and all sufficiently large $M$ the density of $n$
with
 \be \label{highprob}
 \prob \bigg(B_{L_{n,v}}(My,Mv)>M(t-\epsilon)\bigg)>1-\epsilon
 \ee
is at least $1-\epsilon.$
\end{lemma}

\begin{proof}
Fix $\epsilon>0$. If
\be B_{L_{n,v}}(My,M(y+tv))>M(t-\epsilon) \label{stupid1} \ee
and
\be B_{L_{n,v}}(M(v+sw(v)),Mv)>-\epsilon M \label{stupid2} \ee
then
\begin{eqnarray*}
        B_{L_{n,v}}(My,Mv)
 &=&      B_{L_{n,v}}(My,M(v+sw(v)))+B_{L_{n,v}}(M(v+sw(v)),Mv)\\
 &=&      B_{L_{n,v}}(My,M(y+tv))+B_{L_{n,v}}(M(v+sw(v)),Mv)\\
 &>&      M(t-\epsilon)-\epsilon M\\
 &>&      M(t-2\epsilon).
\end{eqnarray*}
The first line follows from Lemma \ref{busemannproperties}, the
second from (\ref{regan}), and the third from (\ref{stupid1}) and
(\ref{stupid2}).

By Lemma \ref{B_n} for any sufficiently large $M$ the density of
$n$ such that (\ref{stupid1}) is satisfied with probability at
least $1-\epsilon$ is at least $1-\epsilon$.
 If $s\neq 0$ then by
Lemma \ref{slope} for any sufficiently large $M$ the density of
$n$ such that (\ref{stupid2}) is satisfied with probability at
least $1-\epsilon$ is at least $1-\epsilon$. If $s=0$ then
$M(v+sw(v))=Mv$ and (\ref{stupid2}) is satisfied for all $M$ and
$n$. As $\epsilon$ is arbitrary the lemma is true.
\end{proof}


%
%

%

\begin{pfof}{Theorem \ref{maingrowth}}
By the definition of $\sides$ for any $k \leq \sides$ we can find
$v_1,\dots v_k$ such that $v_i \in \good$ for all $i$ and the
lines $L_{v_i}$ are all distinct.  The fact that all $v_i \in
\partial \ball$ and that the tangent lines are distinct implies that $t(v_i,v_j)>0$ for
any $i \neq j$. By multiple applications of Lemma
\ref{lem:highprob} there exists $c>0$ such that for all
$\epsilon>0$ there exists $M$ such that the density of $n$ with
 \be \label{l27}
 \prob \left(B_{L_{n,v_i}}(Mv_j,Mv_i)>cM \ \forall \ i \neq j\right)>1-\epsilon.
 \ee
is at least $1-\epsilon$.  We then choose $x_i$ to be the point in
$\z^2$ nearest to $Mv_i$.
 Thus by Lemma \ref{mugconditions} and (\ref{l27}) we have coexistence with
 probability at least $1-\epsilon$.
\end{pfof}

\begin{pfof}{Theorem \ref{mainfpp}}
This follows from Lemma \ref{mugimpliesgeodesics} and Theorem
\ref{maingrowth}.
\end{pfof}

For the following proofs we will use the following notation. For
$(w,z) \in \R^2$ we use the notation $|(w,z)|=\sqrt{w^2+z^2}$ and
$\Ball(c,r)=\{ a \in \R^2: |c-a|<r\}.$ Let $x,v,y \in V$ have
distinct tangent lines. Let $A=A(x,v,y)\subset \partial \ball$ be
the (open) arc of $\partial \ball$ from $x$ to $y$ that contains
$v$. Remember that $G(\0,L_{n,v})$ is the unique geodesic from $\0$ to $L_{n,v}$.
For any $x,v,y \in V$ we consider the event
\be \label{newline1}
G(\0,L_{n,v})\cap \partial (M\ball) \subset M A.
\ee

\begin{lemma} \label{new}
Let $x,v,y \in V$ have distinct tangent lines and let
$\epsilon>0$. There exists $M_0=M_0(\epsilon,x,v,y)$ such that for any $M>M_0$ we have that the
density of $n$ such that
$$\prob(\text{$(\ref{newline1})$ is satisfied})>1-\epsilon$$
is at least $1-\epsilon.$
\end{lemma}

\begin{proof}
First we claim that for fixed $M,n,\epsilon$ and $z \in \partial
(\ball) \setminus A$ that if
\be\tau(0,Mz)\geq \tau(0,Mv))-\epsilon M \label{pilates} \ee
and $Mz \in G(0,L_{n,v})$ then
\begin{eqnarray*}
\tau(0,L_{n,v})
  & = & \tau(0,Mz)+\tau(Mz,L_{n,v})\\
  & = & \tau(0,Mz)+\tau(Mz,L_{n,v})-\tau(Mv,L_{n,v})+\tau(Mv,L_{n,v})\\
  & \geq & \tau(0,Mv) -\epsilon M + B_{L_{n,v}}(Mz,Mv) +\tau(Mv,L_{n,v})\\
  & \geq & \tau(0,Mv) +\tau(Mv,L_{n,v}) + B_{L_{n,v}}(Mz,Mv) -\epsilon M  \\
  & \geq & \tau(0,L_{n,v}) + B_{L_{n,v}}(Mz,Mv) -\epsilon M.
\end{eqnarray*}
The first line is true because $Mz \in G(0,L_{n,v})$.  The third
line is true because of (\ref{pilates}) and the definition of
$B_{L_{n,v}}$. The last line is true because of the subadditivity
of $\tau$. Thus
$$B_{L_{n,v}}(Mz,Mv) \leq \epsilon M.$$

Fix $\{y_i\}_{i=1}^{k}$, $y_i\in \partial( \ball) \setminus A$ for
all $i$, such that for every $z \in \partial( \ball) \setminus A$
there exists $y_i$ with $|y_i-z|<\epsilon/10(T^*(1,0)+T^*(0,1))$.
Next we note that if
 $$B_{L_{n,v}}(My_i,Mv) \geq 10 \epsilon M$$
 and
 $$\tau (My_i,Mz) \leq 2 \epsilon M$$
 then
$$B_{L_{n,v}}(Mz,Mv)=B_{L_{n,v}}(Mz,My_i)+B_{L_{n,v}}(My_i,Mv)
    \geq 10 \epsilon M - \tau (My_i,Mz) > 2\epsilon M.$$
Thus to bound the probability that there exists $z \in \partial(
\ball) \setminus A$ such that $Mz\in G(0,L_{n,v})$ we need only to
bound the probabilities of
\begin{enumerate}
  \item $\tau(0,Mz)\geq \tau(0,Mv))-\epsilon M$ for all $z \in \partial \ball$,
  \item $B_{L_{n,v}}(My_i,Mv) \geq 10 \epsilon M$ for all $y_i$, and
  \item $\tau (My_i,Mz) \leq  2\epsilon M$
  for all $y_i$ and $Mz \in \Ball(My_i,\epsilon M/(T^*(1,0)+T^*(0,1)))$.
\end{enumerate}

For sufficiently large $M$ the first and third events happen with
probability $1-\epsilon$ by Theorem \ref{shape1}.  As
 $y_i \in \partial \ball \setminus A$
we can write $y=(1-t)v+sw(v)$ with $t>0$.
Thus by Lemma \ref{lem:highprob} we have that the  density of $n$ such
that the second event happens with probability at least
$1-\epsilon$ is at least $1-\epsilon$.
\end{proof}

\begin{pfof}{Theorem \ref{mainfpp2}}
Let $v_1,\dots v_k \in \good$ have distinct tangent lines.  By
Lemma \ref{new} we see that for $i=1,\dots,k/2$ there exists $M$
and infinitely many $n$ such that the finite geodesics
$G(\0,L_{n,v_{2i}})$ are pairwise disjoint in $M\partial \ball$.
They all intersect at $\0$ so for infinitely many $n$ the
geodesics are pairwise disjoint in the complement of $M\ball$.  Thus by compactness we can take weak limits to get at least $k$ infinite geodesics that are pairwise disjoint in the complement of $M\ball$.
Thus $|K(\Gamma(0))| \geq k.$
\end{pfof}

To prove Theorem \ref{exponential} we consider the event
\be \label{newline2}
G(Mv,L_{v,n})\cap M\ball \subset \Ball(Mv,\epsilon M).
\ee

\begin{lemma} \label{new2}
Let $\epsilon>0$.  There exists $M_0=M_0(\epsilon)$ such that for any $M>M_0$ we have that
the density of $n$ such that
$$\prob(\text{}(\ref{newline2})\text{ is satisfied})>1-\epsilon$$
is at least $1-\epsilon.$
\end{lemma}

\begin{proof}
If there exists $z \in M\ball \setminus \Ball(Mv,\epsilon M)$ and
$z \in G(Mv,L_{nv})$ then there exists $z \in G(Mv,L_{nv})$ such
that
 $$z \in Z=\partial \bigg(M\ball \setminus \Ball(Mv,\epsilon M)\bigg).$$

Choose $\{y_i\}_{i=1}^k$ such that for any
$z \in Z$ there exists $y_i$ such that
$$|z-y_i|<\epsilon/100(T^*(1,0)+T^*(0,1)).$$
Suppose the following events happen:
\begin{enumerate}
\item $B_{L_{n,v}}(Mv,My_i)< \epsilon M(T^*(1,0)+T^*(0,1))/10$
     for all $i$,\label{one}
\item $\tau (Mv,Mz)>\epsilon M(T^*(1,0)+T^*(0,1))/3$ for all $z$ such that
 $|z -  v| \geq \epsilon $  \label{two} and
\item $\tau(My_i,Mz)< \epsilon M(T^*(1,0)+T^*(0,1))/10$ for all $i$ and $z$ such that
  $|z-y_i|<\epsilon /100(T^*(1,0)+T^*(0,1))$. \label{three}
\end{enumerate}
Then we claim that (\ref{newline2}) is satisfied.  To see this we
note that by \ref{one} and \ref{three}
\begin{eqnarray}
B_{L_{n,v}}(Mv,Mz)
 &=& B_{L_{n,v}}(Mv,My_i)+B_{L_{n,v}}(My_i,Mz) \nonumber\\
 &<& \epsilon M(T^*(1,0)+T^*(0,1))/10+\tau(My_i,Mz) \nonumber\\
 &<& \epsilon M(T^*(1,0)+T^*(0,1))/5. \label{four}
\end{eqnarray}
 for all $z \in Z$.  If $Mz \in G(Mv,L_{n,v})$ then
 $$B_{L_{n,v}}(Mv,Mz)=\tau (Mv,Mz).$$
 Thus by condition \ref{two} if $z \in Z$ and $Mz \in G(Mv,L_{n,v})$
 then
$$B_{L_{n,v}}(Mv,Mz)=\tau (Mv,Mz)>\epsilon M(T^*(1,0)+T^*(0,1))/3$$
which contradicts (\ref{four}) and establishes the claim.

Thus to prove the lemma we need to show that the density of $n$
such that the probability of all of the events in \ref{one},
\ref{two} and \ref{three} occurring is greater than $1-\epsilon.$
By the argument in Lemma \ref{lem:highprob} for sufficiently large
$M$ with probability at least $1-\epsilon/3$ the density of $n$
such that \ref{one} occurs is at least $1-\epsilon$.  By Theorem
\ref{shape1} the probabilities of the last two events can be made
greater that $1-\epsilon/3$.
\end{proof}

For the final proof we will be dealing with multiple realizations
of first passage percolation.  To deal with this we will use the
notation $\tau^{\omega}(x,y)$ $B^{\omega}_S(x,y)$ and
$G^{\omega}(x,y)$ to represent the quantities $\tau(x,y)$
$B_S(x,y)$ and $G(x,y)$ in $\omega$.

\begin{pfof}{Theorem \ref{exponential}}
Given $k$ by Theorem \ref{mainfpp} we can choose $M$ and
$x_1,\dots,x_k \in \partial M\ball$ such that with positive
probability there exist disjoint geodesics $G_i$ starting at each
$x_i$. By Lemma \ref{lem:highprob} we have that there exists a
measurable choice of geodesics $G_i$ and vertices $x_i$ such that
for any $i \neq j$
 \be \label{dumb} \hat B_{G_i}(x_j,x_i)>100.\ee
There exist finite paths $\tilde G_i \subset M\ball$ and an event
$E$ of positive probability that satisfy the following condition.
For each $\omega \in E$ and $i$, the paths $G_i$ and $\tilde G_i$
agree in $M\ball$. Let $y_i \in \Z^2$ be the first vertex in $G_i$
after  $G_i$ exits $M \ball$ for the last time. We can find
$a_i>0$ and restrict to a smaller event of positive probability
where
\be \label{final} a_i < \tau(x_i,y_i)<a_i+1.\ee
We can pick some large $K$ and further restrict our event as
follows.  Let $z,z' \in \Z^2 \setminus M\ball$ be such that there
exist $x,x' \in \Z^2 \cap M\ball$ and $|z-z'|=|x-x'|=1$.  We
require that for any such $z$ and $z'$ that there exists a path
from $z$ to $z'$ that lies entirely outside of $M\ball$ and has
passage time less than or equal to $K$.  For $K$ large enough the
resulting event $\hat E$ will have positive probability. We now
create a new event $E'$ by taking any $\omega \in \hat E$ and
altering the passage times in $M\ball$.  We will do this in a way
such that $E'$ has positive probability and the inequality
$|K(\Gamma(\0))|\geq k$ is satisfied for all $\omega \in E'$.

First we choose paths $\hat G_i \subset M\ball$ that connect $\0$
to $y_i$ such that $\hat G_i \cap \tilde G_j =\emptyset$ for all
$i\neq j$. This is possible by Lemma \ref{new2}.   A configuration
$\omega' \in E'$ if
\begin{enumerate}
\item there exists an $\omega \in \hat E$ such that
$\omega(e)=\omega'(e)$ for all edges $e$ with both endpoints in
$(M\ball)^C$ and
\be
a_i<\tau^{\omega'}(\0,y_i)<a_i+2,\label{dumber}\ee
\item \label{done}
for every $z \in \Z^2 \setminus M\ball$ there exists $i$ such that
$$G(\0,z)|_{M \ball}=\hat G_i,$$
and
\item for every $z \in \Z^2 \setminus M\ball$ and every $i$
$$G(y_i,z) \subset (M \ball)^C.$$
\end{enumerate}
Note that the first and last conditions imply that for every $z
\in \Z^2 \setminus M\ball$ and every $i$
\be \tau^{\omega'}(y_i,z) \geq \tau^{\omega}(y_i,z).\label{ultimate}\ee
Also note that the first condition
implies if  $v \in G_i \setminus M\ball$ then
\be
\tau^{\omega'}(y_i,v) = \tau^{\omega}(y_i,v). \label{last}
\ee

Fix $\omega' \in E'$ and $v \in G_i\setminus M\ball$.  We claim
that $G^{\omega'}(\0,v)=\hat G^i \cup G^{\omega}(y_i,v).$ By
condition \ref{done} we know that $G^{\omega'}(\0,v)$ must pass
through some $y_l$.  Then we calculate
\begin{eqnarray}
\tau^{\omega'}(\0,y_i)+\tau^{\omega'}(y_i,v)
 &<& a_i+2 + \tau^{\omega}(y_i,v)\nonumber\\
 &<& \tau^{\omega}(x_i,y_i)+2 + \tau^{\omega}(y_i,v)\nonumber\\
 &<& \tau^{\omega}(x_i,v)+2.\label{end1}
 \end{eqnarray}
The first inequality follows from (\ref{dumber}) and (\ref{last}).
The second inequality follows from (\ref{final}). The third is
true because $x_i,y_i$ and $v$ are all on $G_i$. Next we calculate
\begin{eqnarray}
\tau^{\omega'}(\0,y_j)+\tau^{\omega'}(y_j,v)
 &>& a_j + \tau^{\omega}(y_j,v)\nonumber\\
 &>& \tau^{\omega}(x_j,y_j) -1+ \tau^{\omega}(y_j,v)\nonumber\\
 &>& \tau^{\omega}(x_j,v)-1\nonumber\\
 &>& \tau^{\omega}(x_i,v)-1+\hat B_{G_i}(x_j,x_i)\nonumber\\
 &>& \tau^{\omega}(x_i,v)+99. \label{end2}
\end{eqnarray}
The first inequality follows from (\ref{dumber}) and
(\ref{ultimate}).  The second inequality follows from
(\ref{final}).  The third from the subadditivity of $\tau$, the
fourth from the definition of $\hat B_{G_i}$ and the final from
(\ref{dumb}).

Combining (\ref{end1}) and (\ref{end2}) we get that for every $v
\in G_i$ the geodesic $G^{\omega'}(\0,v)$ passes through $y_i$. As
this holds true for every $i$ we have that
$$\prob(|K(\Gamma(\0))|\geq k)\geq \prob(E').$$

The conditions on $E'$ can be satisfied by picking the passage
times through the edges in $\cup \hat G_i$ to be in some
appropriate interval to satisfy (\ref{dumber}) and by choosing the
edges not in $\cup \hat G_i$ to have passage times larger than
$10(K+2+\max a_i)$.
As $\mu$ has finite energy we see that
$$\prob(|K(\Gamma(\0))|\geq k)\geq \prob(E')>0.$$
\end{pfof}

\section*{Acknowledgements}
I would like to thank Itai Benjamini for introducing me to this
problem and Oded Schramm for noticing a error in an earlier
version of this paper.

\end{document}